\documentclass[10pt]{article}
\usepackage{microtype}
\usepackage[english,italian]{babel}
\usepackage{amsfonts}
\usepackage{amssymb}
\usepackage{amsmath}
\usepackage{amsthm}
\usepackage{MnSymbol}
\usepackage{enumerate}
\usepackage{leftidx}
\usepackage{geometry}
\usepackage[T1]{fontenc}

\newcommand{\G}{\Gamma}
\renewcommand\mod{\operatorname{mod}}
\newcommand{\ZZ}{\mathbb{Z}}
\newcommand{\cB}{\mathcal B}

\newcommand{\cP}{\mathcal P}

\newcommand{\cL}{\mathcal L}

\title{Down-linking $(K_v,\Gamma)$-designs to $P_3$-designs}

\author{
\begin{otherlanguage*}{italian}
A. Benini, L. Giuzzi and A. Pasotti\thanks{
anna.benini@ing.unibs.it, luca.giuzzi@ing.unibs.it,
anita.pasotti@ing.unibs.it,
Dipartimento di Ma\-te\-ma\-ti\-ca,
Facolt\`a di Ingegneria,
Universit\`a degli Studi di Brescia,
Via Valotti 9,
I-25133 Brescia (IT).
\vskip.2pt
The present research was performed within the activity of GNSAGA
of the Italian INDAM with the financial support of the Italian
Ministry MIUR, projects ``Strutture di incidenza e combinatorie''
 and ``Disegni combinatorici, grafi e loro applicazioni''.
}\end{otherlanguage*}}

\date{}
\newtheorem{defi}{Definition}[section]
\newtheorem{prop}[defi]{Proposition}
\newtheorem{lem}[defi]{Lemma}
\newtheorem{rem}[defi]{Remark}
\newtheorem{ex}[defi]{Example}
\newtheorem{thm}[defi]{Theorem}

\begin{document}
\selectlanguage{english}
\maketitle
\selectlanguage{english}
\begin{abstract}
Let $\Gamma'$ be a subgraph of a graph $\Gamma$.
We define a \emph{down-link} from a $(K_v,\Gamma)$-design $\cB$ to a
$(K_n,\Gamma')$-design $\cB'$ as a map $f:\cB\to \cB'$
mapping any block of $\cB$ into one of its subgraphs.
This is a new concept, closely
related  with both the notion of
\emph{metamorphosis} and that of  \emph{embedding}.
In the present paper we study down-links in general and prove that any
$(K_v,\Gamma)$-design might be down-linked to a $(K_n,\Gamma')$-design,
provided that $n$ is admissible and large enough.
We also show that 
if $\Gamma'=P_3$, it is always possible to find a down-link
to a design of order at most $v+3$.
This bound is then improved
for several classes of graphs $\Gamma$, by providing explicit
constructions.
\end{abstract}

\noindent {\bf Keywords:} down-link;  metamorphosis; embedding;
           $(K_v,\Gamma)$-design.
\par\noindent {\bf MSC(2010):}  05C51, 05B30, 05C38.

\section{Introduction}
Let $K$ be a graph and $\Gamma\leq K$.
A $(K,\Gamma)$-\emph{design}, also called a \emph{$\Gamma$-decomposition} of $K$,
is a set $\cB$ of graphs all isomorphic to $\Gamma$,
called \emph{blocks}, partitioning the edge-set
of $K$.
Given a graph $\Gamma$, the problem of determining the existence
of  $(K_v,\Gamma)$-designs, also called \emph{$\Gamma$-designs of
order $v$}, where $K_v$ is the complete graph on $v$
vertices, 
has been extensively studied;
for surveys on this topic see, for instance, \cite{BO,BZ}.

We propose the following new definition.
\begin{defi}
\label{d1}
  Given a $(K,\Gamma)$-design $\cB$ and a $(K',\Gamma')$-design $\cB'$
  with $\Gamma'\leq \Gamma$, a \emph{down-link} from $\cB$ to $\cB'$
  is a function $f: \cB \rightarrow \cB'$
  such that $f(B)\leq B$, for any $B \in \cB$.
\end{defi}
By the definition of $(K,\Gamma)$-design,
a down-link is necessarily injective.
When a function $f$ as in Definition \ref{d1}  exists,
that is if each block of $\cB$ contains at least one element of $\cB'$
as a subgraph, it will be said that it is possible to
\emph{down-link} $\cB$ to $\cB'$.
\par
In this paper we shall investigate the  existence and some
further properties of down-links between designs on  complete
graphs and
outline their relationship with some previously known notions.
More in detail, 
Section \ref{all} is dedicated to the close interrelationship
between  down-links, metamorphoses and  embeddings.
In Section \ref{spectrum} we will introduce, in close analogy to embeddings,
two  problems on the spectra of down-links and determine bounds
on their minima.
In Section \ref{sec:gen}  down-links from any  $(K_v,\Gamma)$-design
to a $P_3$-design of order $n\leq v+3$ are constructed; this will
improve on the values determined in Section \ref{spectrum}.
In further Sections \ref{star}, \ref{kite}, \ref{cycle}, \ref{path}
 the existence of
 down-links to $P_3$-designs from,
 respectively, star-designs, kite-designs,  cycle systems and
 path-designs are investigated by providing explicit constructions.
\par
Throughout this paper the following standard notations will be used;
see also \cite{Ha}.
For any graph $\Gamma$, write $V(\Gamma)$ for the set of its vertices and
$E(\Gamma)$ for the set of its edges.
By $t\Gamma$ we shall denote the disjoint union of $t$ copies of graphs
all isomorphic to $\Gamma$.
Given any set $V$,  the complete graph with vertex-set $V$ is $K_V$.
As usual,  $K_{v_1,v_2,\ldots,v_m}$ is the complete $m$--partite graph
with parts of size respectively $v_1,\ldots,v_m$; when $v=v_1=v_2=\cdots=v_m$
we shall simply write $K_{m\times v}$.
When we want to focus our attention on the actual parts
$V_1,V_2,\ldots,V_m$, the
notation $K_{V_1,V_2,\ldots,V_m}$ shall be used instead.
The join $\Gamma+\Gamma'$ of two graphs
consists of the graph $\Gamma\cup\Gamma'$ together with  the
edges connecting all the vertices of $\Gamma$ with all the
vertices of $\Gamma'$; hence,
$\Gamma+\Gamma'=\Gamma\cup\Gamma'\cup K_{V(\Gamma),V(\Gamma')}$.

\section{Down-links, metamorphoses, embeddings}
\label{all}
As it will be shown,
the concepts of down-link, metamorphosis and embedding are closely related.
\par
Metamorphoses of designs
have been first introduced by Lindner
and Rosa in \cite{LR1} in the case $\Gamma=K_4$ and
$\Gamma'=K_3$.
In recent years metamorphoses and their generalizations have been extensively
studied; see for instance \cite{CLFT2008,CLFT2007,KL,LR2,LSMeta,LMQ}.
We here recall the general notion of metamorphosis.
\par\noindent
Suppose $\Gamma'\leq \Gamma$ and
let $\cB$ be a
$(\leftidx{^\lambda}K_v,\Gamma)$-design. For each block $B\in \cB$
take a subgraph $B'\leq B$
isomorphic to $\Gamma'$ and put it into a set $S$.
If it is possible
to reassemble all the remaining edges of $\leftidx{^\lambda}K_v$
into a set $R$ of copies of $\Gamma'$, then $S\cup R$ are the blocks
of a $(\leftidx{^\lambda}K_v,\Gamma')$-design,
which is said to be a metamorphosis of $\cB$.
Thus, if $\cB'$ is a metamorphosis of $\cB$ with $\lambda=1$,
then there exists a
down-link $f:\cB\to\cB'$ given by $f(B)=B'$. With a
slight abuse of notation we shall call \emph{metamorphoses}
all down-links from a
$(K,\Gamma)$-design to a $(K,\Gamma')$-design.
\par
There is also a
generalization of metamorphosis, originally from
\cite{LMQ}, which turns out to be closely related to down-links.
Suppose $\Gamma'\leq \Gamma$ and
let $\cB$ be a
$(\leftidx{^\lambda}K_v,\Gamma)$-design. Write $n$ for the minimum integer
$n\geq v$ for which there exists a
$(\leftidx{^\lambda}K_n,\Gamma')$-design.
Take $X=V(K_v)$ and $X\cupdot Y=V(K_n)$.
For each block $B\in \cB$
extract a subgraph of $B$ isomorphic to $\Gamma'$ and put it into a set $S$.
Let also $R$ be the set of all the remaining edges of $\leftidx{^\lambda}K_v$.
Let $T$ be the set of edges of $\leftidx{^\lambda}K_Y$
and of the $\lambda$-fold  complete bipartite graph $\leftidx{^\lambda}K_{X,Y}$.
If it is possible
to reassemble the edges of $R\cup T$ into a set $R'$ of copies of $\Gamma'$,
then $S\cup R'$ are the blocks
of a $(\leftidx{^\lambda}K_n,\Gamma')$-design $\cB'$. In this case,
one speaks of a metamorphosis of $\cB$ into a minimum $\Gamma'$-design.
It is easy to see that for $\lambda=1$ these generalized
metamorphoses also induce down-links.

Even if  metamorphoses with $\lambda=1$  are all   down-links,
the converse is not true.
For instance, all down-links from designs of order $v$ to designs
of order $n<v$ are not metamorphoses. Example \ref{ex:bijective}
shows that such down-links may exist.
\par
Gluing of metamorphoses and down-links can be used to produce new classes of
down-links from old, as shown by the following construction.
Take $\cB$ as a $(K_v,\Gamma)$-design with
$V(K_v)=X\cupdot\bigcupdot_{i=1}^{t}A_i$ and 
suppose $X'\subseteq X$. Let $\G'\leq \G$ and
$\cB'$ be a $(K_{v-|X'|},\G')$-design with $V(K_{v-|X'|})=V(K_{v})\setminus X'$.
Suppose that
\begin{align*}
f_i: (K_{A_i},\G)\textrm{-design}\ \longrightarrow
\ (K_{A_i},\G')\textrm{-design\quad for\ any}\ i=1,\ldots,t,
\end{align*}
\begin{align*}
h_{ij}: (K_{A_i,A_j},\G)\textrm{-design}\ \longrightarrow
\ (K_{A_i,A_j},\G')\textrm{-design\quad for}\ 1\leq i<j\leq t
\end{align*}
are metamorphoses and that
\begin{align*}
g: (K_X,\G)\textrm{-design} \ \longrightarrow 
\ (K_{X\setminus X'},\G')\textrm{-design},
\end{align*}
\begin{align*}
g_{i}: (K_{X,A_i},\G)\textrm{-design}\ 
\longrightarrow \ (K_{X\setminus X',A_i},\G')\textrm{-design\quad for\ any}\
i=1,\ldots, t
\end{align*}
are down-links.
As
$$
K_v=\bigcupdot_{i=1}^t K_{A_i}\cupdot \bigcupdot_{1\leq i<j\leq t}K_{A_i,A_j}
\cupdot K_X 
\cupdot \bigcupdot_{i=1}^t K_{X,A_i}$$
and
$$K_{v-|X'|}=\bigcupdot_{i=1}^t
K_{A_i}\cupdot \bigcupdot_{1\leq i<j\leq t}K_{A_i,A_j}
\cupdot K_{X\setminus X'}
\cupdot
\bigcupdot_{i=1}^t K_{X\setminus X',A_i},$$
the function obtained
by gluing together $g$ and all of the  $f_i$'s, $h_{ij}$'s and $g_{i}$'s 
provides
a down-link from $\cB$ to $\cB'$.

Recall that
an \emph{embedding} of a design $\cB'$
into a design $\cB$ is a function $\psi:\cB' \rightarrow \cB$
such that $\Gamma\leq \psi(\Gamma)$, for any $\Gamma \in \cB'$;
see \cite{Q2002a}.
Existence of embeddings of designs has been widely investigated.
In particular,
a great deal of results are known on injective
embeddings of path-designs;
see, for instance,
\cite{CLQ2005,GQ,MQ99,Q2002,Q2003}.
If $\psi:\cB'\to\cB$ is a \emph{bijective} embedding,
then $\psi^{-1}$ is a down-link from $\cB$ to $\cB'$.
Clearly, a bijective embedding of $\cB'$ into $\cB$ might exist only if
$\cB$ and $\cB'$ have the same number of blocks. This condition, while
quite restrictive, does not necessarily lead to
trivial embeddings, as shown in the following example.
\begin{ex}\label{ex:bijective}
Consider the $(K_4,P_3)$-design
$$\cB'=\{\Gamma'_1=[1,2,3],\Gamma'_2=[1,3,0],\Gamma'_3=[2,0,1]\}$$
 and the $(K_6,P_6)$-design
$$\cB=\{\Gamma_1=[4,0,5,1,2,3],\Gamma_2=[2,5,4,1,3,0],\Gamma_3=[5,3,4,2,0,1]\}.$$
Define $\psi:\cB' \rightarrow \cB$ by
$\psi(\Gamma'_i)=\Gamma_i$ for $i=1,2,3$.
Then,  $\psi$ is a bijective embedding;
consequently, $\psi^{-1}$ is a down-link from $\cB$ to $\cB'$.
\end{ex}

\section{Spectrum problems}\label{spectrum}
Spectrum problems about the existence
of embeddings of designs have been  widely
investigated; see \cite{CLQ2005,DW,GQ,MQ99,Q2002,Q2003}.

In close analogy, we pose  the following questions
about the existence of down-links:
\begin{enumerate}[(I)]
\item  For each admissible $v$, determine the set
  ${\cL}_1\Gamma(v)$ of all integers $n$ such that
  there exists \emph{some} $\Gamma$-design of order $v$ down-linked
  to a $\Gamma'$-design of order $n$.
\item For each admissible $v$, determine the set ${\cL}_2\Gamma(v)$ of all
  integers $n$ such that
  \emph{every} $\Gamma$-design of order $v$ can be down-linked
  to a $\Gamma'$-design
  of order $n$.
\end{enumerate}

In general, write $\eta_i(v;\Gamma,\Gamma')=\inf\cL_i\Gamma(v)$.
When the graphs $\Gamma$ and $\Gamma'$ are easily understood from the
context, we shall simply use $\eta_i(v)$ instead of
$\eta_i(v;\Gamma,\Gamma')$.

The problem of the actual existence of down-links
for given $\Gamma'\leq\Gamma$ is addressed in
Proposition \ref{prop:nonempty}. We recall the following lemma on
the existence of finite embeddings for partial decompositions,
a straightforward consequence of an asymptotic result
by R.M. Wilson \cite[Lemma 6.1]{W}; see also
\cite{BP}.

\begin{lem}
\label{lem:huge}
 Any partial $(K_v,\Gamma)$-design  can be embedded into
  a $(K_n,\Gamma)$-design with
  $n=O((v^2/2)^{v^2})$.
\end{lem}

\begin{prop}
\label{prop:nonempty}
  For any $v$ such that there exists a
  $(K_v,\Gamma)$-design  and any $\Gamma'\leq\Gamma$,
  the sets $\cL_1\Gamma(v)$ and $\cL_2\Gamma(v)$ are non-empty.
\end{prop}
\begin{proof}
Fix first a $(K_v,\Gamma)$-design $\cB$.
Denote by $K_v(\Gamma')$ the so called \emph{complete $(K_v,\Gamma')$-design},
that is the
set of all subgraphs of $K_v$ isomorphic to $\Gamma'$, and let
$\zeta:\cB\to K_v(\Gamma')$ be any function such that $\zeta(\Gamma)\leq\Gamma$
for all $\Gamma\in\cB$. Clearly, the image of $\zeta$ is a partial
$(K_v,\Gamma')$-design $\cP$; see \cite{CHLR}. By Lemma \ref{lem:huge},
there is an integer $n$
such that $\cP$ is  embedded into a $(K_n,\Gamma')$-design $\cB'$.
Let $\psi:\cP\to\cB'$ be such an embedding;
then, $\xi=\psi\zeta$ is, clearly,
a down-link from $\cB$ to a $\Gamma'$-design $\cB'$ of
order $n$. Thus, we have shown that
for any $v$ such that a $\Gamma$-design of order $v$ exists,
and for any $\Gamma'\leq\Gamma$ the set
$\cL_1\Gamma(v)$ is non-empty.
\par
To show that $\cL_2\Gamma(v)$ is also non-empty, proceed as follows.
Let $\omega$ be the number of distinct $(K_v,\Gamma)$-designs $\cB_i$.
For any $i=0,\ldots,\omega-1$, write $V(\cB_i)=\{0,\ldots,v-1\}+i\cdot v$.
Consider now  $\Omega=\bigcup_{i=0}^{\omega-1}\cB_i$.
Clearly, $\Omega$ is a partial $\Gamma$-design of order $v\omega$.
As above, take $K_{v\omega}(\Gamma')$ and construct a function
$\zeta:\Omega\to K_{v\omega}(\Gamma')$ associating to each $\Gamma\in\cB_i$
a $\zeta(\Gamma)\leq\Gamma$. The image $\bigcup_i\zeta(\cB_i)$ is
a partial $\Gamma'$-design
$\Omega'$. Using Lemma
\ref{lem:huge} once more, we determine an integer $n$ and an embedding
$\psi$ of $\Omega'$ into a $(K_n,\Gamma')$-design $\cB'$.
For any $i$, let $\zeta_i$ be the restriction of $\zeta$ to $\cB_i$.
It is straightforward to see that $\psi\zeta_i:\cB_i\to\cB'$
is a down-link from $\cB_i$ to a $(K_n,\Gamma')$-design. It
follows that $n\in\cL_2\Gamma(v)$.
\end{proof}
Notice that the order of magnitude of $n$ is $v^{2v^2}$; yet,
it will be shown
that in several cases it is possible to construct
down-links from  $(K_v,\Gamma)$-designs to $(K_n,\Gamma')$-designs with $n\approx v$.

Lower bounds on $\eta(v;\Gamma,\Gamma')$ are usually hard to obtain
and might not be strict;
a easy one to prove is the following:
\[
(v-1)\sqrt{\frac{|E(\Gamma')|}{|E(\Gamma)|}}<
\eta_1(v;\Gamma,\Gamma'). \]

\section{Down-linking $\Gamma$-designs to $P_3$-designs}
\label{sec:gen}
From this section onwards we shall fix $\Gamma'=P_3$
and focus our attention on the
existence of down-links to $(K_n,P_3)$-designs.
Recall that a $(K_n,P_3)$-design
exists if, and only if, $n\equiv 0,1(\mod 4)$; see \cite{T}.
We shall make extensive use of the following result from
\cite{We}.
\begin{thm}\label{decP3} Let $\Gamma$ be a connected graph. Then, the
edges of $\Gamma$ can be partitioned into copies of $P_3$ if and only if
the number of edges is even. When the number of edges is odd,
$E(\Gamma)$ can be partitioned into a single edge together with copies
of $P_3$.
\end{thm}
Our main result for down-links from a general $(K_v,\Gamma)$-design
is contained in the following theorem. 
\begin{thm}
\label{t:v+3}
  For any $(K_v,\Gamma)$-design $\cB$ with $P_3\leq \Gamma$,
  \[\eta_1(v)\leq \eta_2(v)\leq v+3.\]
\end{thm}
\begin{proof}
For any block $B\in\cB$, fix a $P_3\leq B$ to be used for
the down-link. Write $S$ for the
set of all these $P_3$'s. Remove  the
edges covered by $S$ from $K_v$ and consider the
remaining graph $R$.
If each connected component of $R$  has an even number of edges, by
Theorem \ref{decP3}, there is a decomposition $D$ of $R$ in $P_3$'s;
$S\cup D$ is a decomposition of $K_v$;  thus, $\eta_1(v)\leq\eta_2(v)\leq v$.
If not, take $1\leq w\leq 3$
such that $v+w\equiv 0,1\pmod{4}$. Then, the graph $R'=(K_v+K_w)\setminus S$ 
is connected and  has an even number of edges.
Thus, by Theorem \ref{decP3}, there
is a decomposition $D$ of $R'$  into copies of $P_3$'s. It follows that
$S\cup D$ is a $(K_{v+w},P_3)$-design $\cB'$.
\end{proof}
\begin{rem}
In Theorem \ref{t:v+3}, if $v\equiv 2,3\pmod{4}$, then
the order of the design $\cB'$ is the 
smallest $m\geq v$ for which there exists a $(K_m,P_3)$-design.
Thus, 
the down-links are actually metamorphoses to minimum $P_3$-designs.
This is not the case for $v\equiv 0,1\pmod{4}$, as we cannot 
\emph{a priori}
guarantee 
that each connected component of $R$ has an even
number of edges.
\end{rem}
Theorem \ref{t:v+3} might be improved under some further (mild)~assumptions
on $\Gamma$.
\begin{thm}
\label{t:44}
Let $\cB$ be a $(K_v,\Gamma)$-design.
\begin{enumerate}[a)]
\item
If $v\equiv 1,2\pmod{4}$, $|V(\Gamma)|\geq 5$ and there are at least $3$
vertices in $\Gamma$ with degree at least $4$, then there exists a
down-link from $\cB$ to a $(K_{v-1},P_3)$-design.
\item
If $v\equiv 0,3\pmod{4}$, $|V(\Gamma)|\geq 7$ and there are at least $5$
vertices in $\Gamma$ with degree at least $6$, then there exists a
down-link from $\cB$ to a $(K_{v-3},P_3)$-design.
\end{enumerate}
\end{thm}
\begin{proof}
\begin{enumerate}[a)]
\item\label{dd:a}
Let $x,y\in V(K_v)$. Extract from any $B\in\cB$
a $P_3\leq B$ whose
vertices are neither $x$ nor $y$ and use it for the down-link.
This is always possible,
since $|V(\Gamma)|\geq 5$ and there is at least one vertex in 
$\Gamma\setminus\{x,y\}$ of degree at least $2$. Write now $S$
for the set of all of these $P_3$'s. 
 Consider the graph $R=(K_{v-2}+\{\alpha\})\setminus S$ where
 $K_{v-2}=K_{v}\setminus\{x,y\}$.
 This is a connected graph with an even number of edges; thus,
 by Theorem \ref{decP3}, there exists a decomposition $D$ of $R$
 in $P_3$'s. Hence, $S\cup D$ provides the blocks of a
 $P_3$-design of order $v-1$.
\item
  In this case
  consider $4$ vertices $\Lambda=\{x,y,z,t\}$ of $V(K_v)$. By the assumptions,
  it is always possible to take a $P_3$  disjoint from $\Lambda$
  from each block of $\cB$. We now argue as in the proof of part \ref{dd:a}).
\end{enumerate}
\vskip-.4cm
\end{proof}

The down-links constructed above are not, in general,
to designs whose order is as small as possible;
thus, theorems \ref{t:v+3} and \ref{t:44} do not provide the exact value of
$\eta_1(v)$, unless further assumptions are made.
\begin{rem}
\label{r:emb}
  In general, a $(K_n,P_3)$-design can be trivially
  embedded into  $P_3$-designs of any admissible order $m\geq n$.
  Thus, if $n\in\cL_i\Gamma(v)$, then
  $\{ m\geq n\ | \ m\equiv 0,1 (\mod 4) \} \subseteq \cL_i\Gamma(v)$.
  Hence,
  \[ \cL_i\Gamma(v)=\{ m\geq\eta_i(v)\ | \ m\equiv 0,1 (\mod 4)\}. \]
  Thus,
  solving problems {\rm (I)} and {\rm (II)}
  turns out to be  actually equivalent
  to  determining exactly
  the values of $\eta_1(v;\Gamma,P_3)$ and $\eta_2(v;\Gamma,P_3)$.
\end{rem}
For the remainder of this paper, we shall always silently apply
Remark \ref{r:emb} in all the proofs.

\section{Star-designs}
\label{star}
In this section the existence of down-links from star-designs to
$P_3$-designs is investigated.
We follow  the notation introduced in Section \ref{spectrum},
where $\Gamma'=P_3$ is understood.
 Recall that the \emph{star} on $k+1$
vertices $S_k$ is the complete bipartite graph $K_{1,k}$
with one part having a single vertex, say $c$, called the \emph{center} of
the star, and the other part having $k$ vertices, say $x_i$ for
$i=0,\ldots,k-1$, called \emph{external vertices}.
In general, we shall write $S_k=[c;x_0,x_1,\ldots,x_{k-1}]$.
\par
In \cite{Tstar}, Tarsi proved that a $(K_v,S_k)$-design exists if, and only if,
$v\geq 2k$ and $v(v-1)\equiv 0\pmod{2k}$.
When
$v$ satisfies these necessary conditions
we shall determine the sets $\cL_1S_k(v)$ and
$\cL_2S_k(v)$.

\begin{prop}
\label{r:star}
For any admissible $v$ and $k>3$,
\begin{eqnarray}
  \cL_1S_k(v)&\subseteq& \{n\geq v-1\,|\,
  n\equiv 0,1\ (\mod{4}) \}, \label{e:S1} \\
  \cL_2S_k(2k)&\subseteq&\{n\geq v-1\,|\,
  n\equiv 0,1\ (\mod{4}) \}, \label{e:S2} \\
  \cL_2S_k(v)&\subseteq&\{n\geq v\,|\,
  n\equiv 0,1\ (\mod{4}) \} \text{\ for $v>2k$. } \label{e:S3}
\end{eqnarray}
\end{prop}
\begin{proof}
  In a $(K_v,S_k)$-design $\cB$,
  the edge $[x_1,x_2]$ of $K_v$ belongs either to a star of
  center $x_1$ or to a star of center $x_2$. Thus, there is
  possibly at most one vertex which is not the center of any star;
  \eqref{e:S1} and \eqref{e:S2} follow.
  \par
  The condition \eqref{e:S3} is obvious when
  any vertex of $K_v$ is center of at least one star of $\cB$.
  Suppose now that there exists a vertex, say $x$, which is not center
  of any star. Since $v>2k$, there exists also a vertex $y$ which is
  center of at least two stars. Let $S=[y;x,a_1,\ldots,a_{k-1}]$
  and take, for any $i=1,\ldots,k-1$,
  $S^i$ as the star with center $a_i$ and containing $x$. 
  Replace $S$ in $\cB$ with the star
  $S'=[x;y,a_1,\ldots,a_{k-1}]$.  Also, in each $S^i$ substitute the
  edge $[a_i,x]$ with $[a_i,y]$. Thus, we have again a
  $(K_v,S_k)$-design in which each vertex of $K_v$ is the
  center of at least
  one star. This gives \eqref{e:S3}.
\end{proof}

\begin{thm}\label{thm:star1}
Assume $k>3$. For every $v\geq 4k$ with $v(v-1)\equiv 0(\mod 2k)$,
\[\cL_1S_k(v)=
 \{n\geq v-1\ |\ n\equiv 0,1\ (\mod{4})\}. \]
\end{thm}
\begin{proof}
By Proposition \ref{r:star}, it is enough to show
$\{n\geq v-1\ |\ n\equiv 0,1\pmod{4}\}\subseteq\cL_1S_k(v)$.
We distinguish some cases:
\begin{enumerate}[a)]
\item\label{s:a} $v\equiv 0\pmod{4}$. Since $v$ is admissible and $v\geq 4k$,
by \cite[Theorem 1]{C} there always exists a $(K_v,S_k)$-design
$\cB$ having exactly one vertex, say $x$, which is not the center
of any star. Select from each block of $\cB$ a path $P_3$ whose vertices are
different from $x$. Use these $P_3$'s for the down-link and remove their
edges from $K_v$. This yields a connected graph $R$ having an even
number of edges. So, by Theorem \ref{decP3} $R$ can be decomposed in
$P_3$'s; hence, there exists a down-link from $\cB$ to a
$(K_v,P_3)$-design.
\item\label{s:b} $v\equiv 1,2\pmod{4}$. In this case
  there always exists a $(K_v,S_k)$-design $\cB$ having exactly one vertex,
  say $x$, which is not center of any star and at least one vertex $y$
  which is center of exactly one star, say $S$; see \cite[Theorem 1]{C}. 
  Choose  a $P_3$, say
  $P=[x_1,y,x_2]$, in $S$.
  Let now $S'$ be the star containing the edge $[x_1,x_2]$ and
  pick  a $P_3$ containing this edge. 
  Select from each of the other blocks of $\cB$ a $P_3$
  whose vertices are different from $x$ and $y$.
  This is always possible since $k>3$.
  Use all of these $P_3$'s to construct a down-link.
  Remove from $K_v\setminus\{x\}$ all of the edges of the $P_3$'s,
  thus obtaining a graph $R$ with an even number of edges.
  Observe that $R$ is connected, as $y$ is adjacent to all
  vertices of $K_v$ different from $x,x_1,x_2$.
  Thus, by Theorem \ref{decP3}, $R$ can be decomposed
  in $P_3$'s. Hence, there exists a down-link from
  $\cB$ to a $(K_{v-1},P_3)$-design.
\item\label{s:c} $v\equiv 3\pmod{4}$. As neither $n=v-1$ nor $n=v$ are
admissible for $P_3$-designs,  the result follows 
arguing as in the proof of Theorem
\ref{t:v+3}.
\end{enumerate}
\vskip-1.5\baselineskip
\end{proof}

The condition $v\geq 4k$ might be relaxed 
when $k>3$ is a prime power, as shown by the following theorem.
\begin{thm}  Let $k>3$  be a prime power.
 For every $2k\leq v<4k$ with $v(v-1)\equiv 0\pmod{2k}$,
 \[\cL_1S_k(v)= \{n\geq v-1\ |\ {n\equiv{0,1}\ (\mod{4})}\}. \]
\end{thm}
\begin{proof} Since $k$ is a prime power, $v$ can only assume 
  the following values: $2k, 2k+1, 3k, 3k+1$.  For each of the
  allowed   values of $v$ there exists a $(K_v,S_k)$-design with
  exactly one vertex  which is not center of any
  star; see \cite{C}. The result can be obtained arguing as in
  previous theorem.
\end{proof}

\begin{thm}\label{thm:star}
Let $k>3$ and take $v$ be such that $v(v-1)\equiv 0\pmod{2k}$. Then,
 \[ \cL_2S_k(2k) =
 \{n\geq v-1\ |\ n\equiv 0,1\ (\mod{4})\};\]
 \[ \cL_2S_k(v) =
 \{n\geq v\ |\ n\equiv 0,1\ (\mod{4})\}\ \text{for $v>2k$.} \]
\end{thm}
\begin{proof}
Let $\cB$ be a $(K_{2k},S_k)$-design.
Clearly, there is exactly one vertex
of $K_{2k}$ which is not the center of any star.
By Proposition \ref{r:star}, it is enough to show that
$\{n\geq v-1\ |\ n\equiv 0,1\pmod{4}\}\subseteq\cL_2S_k(2k)$.
The result can be obtained arguing
as in step \ref{s:a}) of Theorem \ref{thm:star1}
for $k$ even and as in step \ref{s:b}) of the same for $k$ odd.
\par
We now consider the case $v>2k$. 
As before,
by Proposition \ref{r:star}, we just need to prove one of the inclusions.
Suppose $v\equiv0,1\pmod{4}$.
Let $\cB$ be a $(K_v,S_k)$-design.
For $k$ even, each star is a
disjoint union of $P_3$'s and  the existence of a down-link to
a $(K_v,P_3)$-design is trivial.
For $k$ odd, observe 
that $\cB$ contains
an even number of stars. Hence,
there is an even number of vertices
$x_0,x_1,x_2,\ldots,x_{2t-1}$ of
$K_v$ which are center of an odd number of stars. 
Consider the edges $[x_{2i},x_{2i+1}]$ for $i=0,\ldots,t-1$.
From each star of $\cB$, extract a $P_3$ which does not contain
any of the aforementioned edges and use it for the down-link. 
If $y\in K_v$ is the center of an even number of
stars, then the union of all the remaining edges of stars with
center $y$ is a connected graph with an even number of edges;
thus, it is possible to apply Theorem \ref{decP3}.
If $y$ is the center of an odd number of stars, then there is
an edge $[x_{2i},x_{2i+1}]$ containing $y$. In this case the
graph obtained by
the union of all the remaining edges of the
stars with centers $x_{2i}$ and $x_{2i+1}$ is connected and has
an even number of edges. Thus, we can apply again Theorem \ref{decP3}.
For $v\equiv 2,3\pmod{4}$, the result follows as in Theorem \ref{t:v+3}.
\end{proof}

\section{Kite-designs}
\label{kite}
Denote by
$D=[a,b,c \bowtie  d]$ the \emph{kite}, a triangle with
an attached edge,
having vertices $\{a,b,c,d\}$ and edges
$[c,a],[c,b],[c,d],$ $[a,b]$.

In \cite{BS}, Bermond and Sch\"{o}nheim
proved that a kite-design of order $v$ exists if, and only if,
 $v\equiv 0,1 (\mod 8)$, $v>1$.
In this section we completely determine the sets $\cL_1D(v)$ and  $\cL_2D(v)$
where $\Gamma'=P_3$ and $v$, clearly, fulfills the aforementioned condition.

We need now to recall some
preliminaries on difference families.
For general definitions and in depth discussion, see \cite{BPBica}.
Let $(G,+)$ be a group and take $H\leq G$.
A set $\cal F$ of kites with vertices in
$G$ is called a
$(G,H,D,1)$-\emph{difference family}
(DF, for short), if the list
$\Delta{\cal F}$ of differences from $\cal F$, namely the list of all
possible differences $x-y$, where $(x,y)$ is an ordered pair of
adjacent vertices of a kite in $\cal F$, covers all the elements of
$G \setminus H$ exactly once,
while no element of $H$ appears in $\Delta{\cal F}$.
\begin{prop}\label{prop:kite}
For every $v\equiv 0,1(\mod\ 8)$, $v>1$,
\begin{eqnarray}
\cL_1D(v) &\subseteq&\{n\geq v-1 |\ n\equiv 0,1\ (\mod{4})\}, \label{l1d}
\\
\cL_2D(v) &\subseteq&\{n\geq v\ |\ n\equiv 0,1\ (\mod{4})\}. \label{l2d}
\end{eqnarray}
\end{prop}
\begin{proof}
Let $\cB$ and $\cB'$ be respectively a $(K_v,D)$-design and 
a $(K_n,P_3)$-design.
Suppose there are $x,y\in V(K_v)\setminus V(K_n)$ with $x\neq y$.
Since 
there is at least one kite $D\in\cB$ containing both $x$ and $y$,
we see that it is not possible to extract any $P_3\in\cB'$ from $D$;
thus $n\geq v-1$. This proves \eqref{l1d}.
\par
As for \eqref{l2d}, we distinguish two cases.
For $v\equiv0\pmod{8}$, there does not exist a
$P_3$-design of order $v-1$.
On the other hand,
 for any $v=8t+1$,
\[ {\cal F}=\{ [2i-1,3t+i,0\bowtie 2i]\, |\, i=1,\ldots,t \} \]
is a $(\ZZ_{8t+1},\{0\},D,1)$-DF. 
As a special case of a more general result proved in \cite{BPBica},
the existence of such a difference family implies that of a cyclic
$(K_{8t+1},D)$-design $\cB$. Thus, any $x\in V(K_{8t+1})$ has degree
$3$ in at least one block of $\cB$. Hence, there is no down-link of
$\cB$ in a design of order less than $8t+1$.
\end{proof}

\begin{lem}\label{Buratti}
For every integer $m=2n+1$ there exists a $(K_{m\times 8},D)$-design.
\end{lem}
\begin{proof}
The set
$$ {\cal F}=\left\{ [(0,0),(0,2i),(2,i) \bowtie (1,0)],
  [(0,0),(4,i),(1,-i) \bowtie (6,i)]\, |\, i=1,\ldots,n \right\}$$
is  a
 $(\ZZ_8 \times \ZZ_m,\ZZ_8 \times \{0\},D,1)$-DF.
A special case of a result in \cite{BPBica} shows that
any difference family with these parameters
determines a $(K_{m\times 8},D)$-design.
\end{proof}

\begin{prop}\label{prop:degree2}
There exists a $(K_v,D)$-design with a vertex $x$ having degree $2$
in all the blocks in which it appears if and only if $v\equiv 1\pmod{8}$, $v>1$.
\end{prop}
\begin{proof}
Clearly, $v\equiv1\pmod{8}$, $v>1$, is a necessary condition for
the existence of such a design. We will show that it is also sufficient.
Assume $v=8t+1$, $t\geq 1$.
Let $A_i=\{a_{i1},a_{i2},\dots,a_{i8}\}$, $i=1,\dots,t$ and write
$V(K_v)=\{0\}\cupdot{A_1}\cupdot{A_2}\cupdot\cdots\cupdot{A_t}$.
Clearly, $E(K_v)$ is the disjoint union
of the sets of edges of $K_{0,A_i}$, $K_{A_i}$ and
$K_{A_1,A_2,\dots,A_t}$, for $i=1,2,\dots,t$.
\begin{list}{\labelitemi}{
\leftmargin=0pt}
\item
  Suppose $t=1$, so that $V(K_v)=\{0\}\cupdot{A_1}$.
  An explicit kite decomposition of $K_v=K_{0,A_1}\cup K_{A_1}$ where
  the degree of $0$ is always 2 is given by
{\small
\begin{align*}
\{ & [0,a_{11},a_{12}\bowtie a_{16}],[0,a_{14},a_{13}\bowtie a_{15} ], [0,a_{16},a_{15}\bowtie a_{17}],
[0,a_{17},a_{18}\bowtie a_{16}], \\ &
[a_{11},a_{13},a_{16}\bowtie a_{17}], [a_{11},a_{17},a_{14}\bowtie a_{16}], [a_{11},a_{15},a_{18}\bowtie a_{14}], \\ &
[a_{12},a_{18},a_{13}\bowtie a_{17}], [a_{14},a_{15},a_{12}\bowtie a_{17}]
\}.
\end{align*}
}
\item Let now $t=2$, so that $V(K_v)=\{0\}\cupdot{A_1}\cupdot{A_2}$.
There exists  a kite decomposition of $K_v$ where the degree of $0$ is
always 2. Such a decomposition
 results from
 the disjoint union of the previous kite decomposition of
 $K_{0,A_1}\cup K_{A_1}$ and
 the kite decomposition of $K_{0,A_2}\cup K_{A_2}\cup K_{A_1,A_2}$ here listed:
{\small
{\begin{align*}
\{ & [0,a_{22},a_{21}\bowtie a_{18}],[0,a_{24},a_{23}\bowtie a_{18}], [0,a_{26},a_{25}\bowtie a_{18}],
[0,a_{28},a_{27}\bowtie a_{18}], \\ & [a_{11},a_{21},a_{28}\bowtie a_{18}], [a_{11},a_{27},a_{22}\bowtie a_{18}],
[a_{11},a_{23},a_{26}\bowtie a_{18}], [a_{11},a_{25},a_{24}\bowtie a_{18}], \\ & [a_{12},a_{27},a_{21}\bowtie a_{17}],
[a_{12},a_{26},a_{22}\bowtie a_{17}], [a_{12},a_{25},a_{23}\bowtie a_{17}], [a_{12},a_{28},a_{24}\bowtie a_{17}],\\ &
[a_{13},a_{21},a_{26}\bowtie a_{17}],[a_{13},a_{22},a_{25}\bowtie a_{17}], [a_{13},a_{23},a_{28}\bowtie a_{17}],
[a_{13},a_{24},a_{27}\bowtie a_{17}], \\ & [a_{14},a_{25},a_{21}\bowtie a_{16}], [a_{14},a_{24},a_{22}\bowtie a_{16}],
[a_{14},a_{23},a_{27}\bowtie a_{16}], [a_{14},a_{26},a_{28}\bowtie a_{16}], \\& [a_{15},a_{24},a_{21}\bowtie a_{23}],
[a_{15},a_{23},a_{22}\bowtie a_{28}], [a_{15},a_{28},a_{25}\bowtie a_{16}], [a_{15},a_{26},a_{27}\bowtie a_{25}],\\ &
[a_{24},a_{26},a_{16}\bowtie a_{23}]
\}.
\end{align*}}}
\item Take $v=8t+1$, $t\geq 3$.
For odd $t$, the complete multipartite graph
$K_{t \times 8}$ always admits a kite decomposition; see Lemma \ref{Buratti}.
Thus, $K_v$ has a kite decomposition which
is the disjoint union of the kite decomposition of
$K_{0,A_i}\cup K_{A_i}$,
for each $i=1,\dots,t$ (compare this with the case $t=1$), and that  of
$K_{A_1,A_2,\dots,A_t}$.
If $t$ is even,
 write
 \[ V(K_v)=\{0\}\cupdot A_1\cupdot\cdots\cupdot{A_{t-1}}\cupdot A_t. \]
As $t-1$ is odd, the graph
$K_v$ has a kite decomposition which is the disjoint
union of the kite decompositions of
$K_{0,A_t}\cup K_{A_t}$ (see the case $t=1$), $K_{A_1,A_2,\dots,A_{t-1}}$
and $K_{0,A_i}\cup K_{A_i}\cup K_{A_i,A_t}$ (as in the case $t=2$),
for  $i=1,\dots,t-1$.
In either case the degree of $0$ is $2$.
\end{list}
\vskip-0.70cm
\end{proof}

\begin{thm}
For every $v\equiv 0,1\pmod{8}$, $v>1$,
 \begin{align}\label{LD1}
 \cL_1D(v) &= \{n\geq v-1\ |\ n\equiv 0,1\ (\mod{4})\};\\
\label{LD2}
 \cL_2D(v) &= \{n\geq v\ |\ n\equiv 0,1\ (\mod{4})\}.
\end{align}
\end{thm}
\begin{proof}
Relation \eqref{LD1} follows from
Proposition \ref{prop:kite} and Proposition \ref{prop:degree2}.
Clearly, any $D$-design with a vertex $x$ of degree $2$
in all the blocks in which it appears can
be down-linked to a $P_3$-design of order $v-1$.
\par
In view of Proposition \ref{prop:kite}, to prove relation \eqref{LD2},
it is sufficient to observe that each kite can be seen as the 
union of two $P_3$'s.
\end{proof}

\section{Cycle systems}\label{cycle}
Denote by $C_k$ the cycle on $k$ vertices, $k\geq 3$.
It is well known that a $k$-cycle system of order $v$, that is
a $(K_v,C_k)$-design, exists if, and only if,
$k\leq v$, $v$ is odd and $v(v-1)\equiv 0\pmod{2k}$.
The \emph{if part} of this theorem was solved by Alspach
and Gavlas \cite{AGA} for $k$ odd and by \~Sajna \cite{SJ}
for  $k$ even.

In this section we shall provide some partial results on
$\cL_1C_k(v)$ and $\cL_2C_k(v)$.

\begin{thm}\label{thm:cycle}
For any admissible $v$,
\begin{align*}
\cL_2C_{3}(v)&=\cL_1C_{3}(v)=\{n\geq v\ |\ n\equiv0,1\ (\mod{4})\};\\
\cL_2C_4(v) &= \cL_1C_4(v)=\{n\geq v-1\ |\ n\equiv 0,1\ (\mod{4})\};\\
\cL_2C_{5}(v)&=\{n\geq v-1\ |\ n\equiv0,1\ (\mod{4})\}\subseteq\cL_1C_{5}(v);
\end{align*}
and for any $k\geq6$
\[\left\{ n\geq v- \left\lfloor\frac{k-4}{3}\right\rfloor\,\big|\, n\equiv0,1\
 (\mod{4})\right\}\subseteq \cL_2C_{k}(v)\subseteq
\cL_1C_{k}(v).\]
\end{thm}
\begin{proof}
\begin{list}{\labelitemi}{
\leftmargin=0pt}
\item Suppose $k=3$.
It is obvious that a down-link from a $(K_v,C_3)$-design $\cB$
to a $P_3$-design
of order less than $v$ cannot exists.
When $v\equiv 1\pmod{4}$, the triangles in $\cB$
can be paired so that each pair share a vertex; see \cite{HR}.
Let $T_1=(1,2,3)$ and $T_2=(1,4,5)$ be  such a pair.
Use the paths $[1,2,3]$ and $[1,4,5]$ for down-link and consider
the path $[3,1,5]$. Observe that these three paths provide a 
decomposition of the edges of $T_1\cup T_2$ in $P_3$'s.
The proof is completed by repeating
this procedure for all paired triangles.
For  $v\equiv 3\pmod{4}$, proceed as in Theorem \ref{t:v+3}.
\item
Assume $k=4$. Let $\cB$ be a $(K_v,C_4)$-design.
It is easy to see that, as in the case of the kites, the image of a $C\in\cB$
in a $(K_n,P_3)$-design $\cB'$ must necessarily leave out exactly one
of the vertices of $C$. Obviously, any two vertices of $V(K_v)$ are
contained together in at least one  block of $\cB$; thus,
$\cL_2C_4(v)\subseteq\cL_1C_4(v)\subseteq \{n\geq v-1\,|\, n\equiv 0,1\pmod{4}\}$.
We now prove the reverse inclusion:
take  $x\in V(K_v)$ and delete from each $C\in\cB$ with $x\in V(C)$
the edges passing through $x$, thus obtaining a $P_3$, say $P$. Let the
image of $C$ under the down-link be $P$.
Observe now that the blocks not containing $x$
can still be decomposed into two $P_3$'s. Thus,  it is possible
to construct a down-link from $\cB$ to
a $P_3$-design of order $v-1$.
\item
Take $k=5$. 
Note that a $(K_v,C_5)$-design $\cB$  necessarily
satisfies either of the following:
\begin{enumerate}[1)]
\item\label{p:1} there exist $x,y\in V(K_v)$ such that $x$ and $y$
appear in exactly one block $B$, wherein they are adjacent;
\item\label{p:2} every pair of vertices of $K_v$ appear in exactly $2$ blocks.
In other words,  $\cB$ is a Steiner pentagon system; see \cite{LS}.
\end{enumerate}
We will show that it is always possible to down-link $\cB$
to a $P_3$-design of order $n\geq v-1$ admissible.
Suppose $v\equiv1\pmod{4}$.
If $\cB$ satisfies \ref{p:1}), then
select from each block a $P_3$ whose vertices
are different from $x$ and $y$.
Use these $P_3$'s to construct the down-link.
Observe that $K_{v-1}=K_v\setminus\{x\}$ minus the edges used for the down-link
is a connected graph; thus the assertion follows from Theorem \ref{decP3}.
If $\cB$  satisfies \ref{p:2}),  
select from each block a $P_3$ whose vertices
are different from $x$ and $y$.
Note that this is always possible, unless the cycle is
$C=(x,a,b,y,c)$. In this case, select from $C$
the path $P=[a,b,y]$. Note that none of the selected paths contains
$x$; thus, their union is a subgraph $S$ of $K_{v-1}=K_v\setminus\{x\}$.
It is easy to see that
each vertex $v\neq b$ of $K_{v-1}\setminus S$ is adjacent
to $y$. Thus, either $K_{v-1}\setminus S$ is connected or it consists
of the isolated vertex $b$ and a connected component.  
In both cases it is
possible to apply Theorem \ref{decP3} to obtain
a $(K_{v-1},P_3)$-design.
When
$v\equiv 3\pmod{4}$, argue as in Theorem \ref{t:v+3}.
\item
Let $k\geq6$ and denote by $\cB$  a $(K_v,C_k)$-design.
Write $t=\left\lfloor\frac{k-4}{3}\right\rfloor$.
Take $t+1$ distinct vertices $x_1,x_2,\ldots,x_t,y\in V(K_v)$. 
Observe that it is always possible to extract from each block $C\in\cB$
a $P_3$ whose vertices are different from $x_1,\ldots,x_t,y$, as
we are forbidding at most $2\lfloor\frac{k-4}{3}\rfloor+2$
edges from any $k$-cycle; consequently,
the remaining edges cannot be pairwise disjoint.
Use these $P_3$'s for the down-link. 
Write $S$ for
the image of the down-link,
regarded as a subgraph of $K_{v-t}=K_v\setminus\{x_1,x_2,\ldots,x_t\}$.
Observe that the edges of $K_{v-t}$ not contained in $S$
form a connected graph $R$.
When $R$ has an even number of edges,
namely $v-\lfloor\frac{k-4}{3}\rfloor\equiv0,1\pmod{4}$, the result is
a direct consequence of Theorem \ref{decP3} and we are done.
Otherwise add $u=1$ or $u=2$ vertices to $K_{v-t}$ and then apply
Theorem \ref{decP3} to $R'=(K_{v-t}+K_u)\setminus S$.
\end{list}
\vskip-\baselineskip
\end{proof}

\begin{rem}
It is not possible to down-link a $(K_v,C_5)$-design with
Property   \ref{p:2}) 
to  $P_3$-designs of order smaller than
$v-1$. On the other hand if a $(K_v,C_5)$-design enjoys Property \ref{p:1}),
then it might be possible to obtain a
down-link to a $P_3$-design of order smaller than $v-1$, as shown
by the following example.
\end{rem}
\begin{ex}
Consider the cyclic $(K_{11},C_5)$-design $\cB$ presented in \cite{BDF}:
{\begin{scriptsize}
\begin{eqnarray*}
\cB&=\{&[0,8,7,3,5],[1,9,8,4,6],[2,10,9,5,7],[3,0,10,6,8], 
[4,1,0,7,9], 
 [5,2,1,8,10], \\ && [6,3,2,9,0],  
[7,4,3,10,1], 
[8,5,4,0,2],[9,6,5,1,3], 
 [10,7,6,2,4]\quad\}.
\end{eqnarray*}
\end{scriptsize}
}
Note that $0$ and $1$ appear together in exactly one block.
It is possible to down-link $\cB$ to the following $P_3$-design of order $9$:
\begin{scriptsize}
\begin{eqnarray*}
\cB'&=\{\hspace{-5pt}&
[8,7,3],[8,4,6],[9,5,7],[10,6,8],[7,9,4],[8,10,5], [6,3,2], 
[4,3,10],[8,5,4], \\ && [3,9,6],  [4,10,7]\ \}\cup 
\{\ [3,5,2], [3,8,9],[7,2,10], 
[10,9,2],[6,7,4],
[4,2,8],[2,6,5]\ \}.
\end{eqnarray*}
\end{scriptsize}
\end{ex}

\section{Path-designs}
\label{path}
In \cite{T}, Tarsi proved that the necessary conditions
for the existence of a $(K_v,P_k)$-design, namely
$v\geq k$ and $v(v-1)\equiv 0\pmod{2(k-1)}$,
are also sufficient.
In this section we investigate
down-links from  path-designs to $P_3$-designs and
provide partial results for $\cL_1P_k(v)$ and $\cL_2P_k(v)$.

\begin{thm}
\label{t:81}
For any admissible $v>1$,
\begin{eqnarray}
\cL_1P_4(v)  & = & \{n\geq v-1\ |\ n\equiv0,1\ (\mod{4})\}; \label{P14}\\
\cL_2P_4(v)  & \subseteq & \{n\geq v\ |\ n\equiv0,1\ (\mod{4})\}.\label{P15}
\end{eqnarray}
\end{thm}
\begin{proof}
Let $\cB$ and $\cB'$ be respectively a $(K_v,P_4)$-design and 
a $(K_n,P_3)$-design.
Suppose there exists a down-link $f:\cB\rightarrow\cB'$.
Clearly, $n>v-2$.
Hence,
$\cL_2P_4(v) \subseteq \cL_1P_4(v)
\subseteq \{n\geq v-1\ |\ n\equiv0,1\pmod{4}\}$.
\par
To show the
reverse inclusion in \eqref{P14} we prove
the actual existence of designs
providing down-links attaining the minimum.
For the case
$v\equiv 1,2\pmod4$ we refer to Subsection \ref{ss:app}.  For $v\equiv3\pmod4$, it
is possible to argue as in Theorem \ref{t:v+3}.
For $v\equiv0\pmod4$, observe that a
$(K_v,P_4)$-design exists if, and only if
$v\equiv0,4\pmod{12}$. 
In particular, for $v=4$, the existence of a down-link from a $(K_4,P_4)$-design to
a $(K_4,P_3)$-design is trivial.
For $v>4$,
arguing as in Subsection \ref{ss:app}, we can obtain a
$(K_v,P_4)$-design $\cB$ with a vertex $0\in V(K_v)$ having 
degree $1$ in each block wherein it appears.
Hence, it is possible to choose for the 
down-link a $P_3$ not containing $0$  from any block of $\cB$.
Denote by $S$ the set of all of these $P_3$'s and consider the
complete graph $K_{v-1}=K_v\setminus\{0\}$. Let now 
$R=(K_{v-1}+\{\alpha\})\setminus S$. Clearly, $R$ is a connected
graph with an even number of edges. Hence, by
Theorem \ref{decP3}, $\eta_1(v)=v$.
\par
In order to prove \eqref{P15},
it is sufficient to show that for any admissible $v$ there exists
a $(K_v,P_4)$-design $\cB$ wherein no vertices can be deleted.
In particular, this is the case if each vertex of $K_v$ has degree $2$
in at least one block of $\cB$.
First of all note that
in a $(K_v,P_4)$-design there is at most one vertex with degree $1$ in each
block where it appears.
Suppose that there actually exists
a $(K_v,P_4)$-design $\overline{\cB}$ with a vertex $x$ as above.
It is not hard to see that there is in $\overline{\cB}$
at least one block 
$P^1=[x,a,b,c]$ 
such that the vertices $a$ and $b$ have degree $2$
in at least another block.
Let $P^2=[x,c,d,e]$.
By reassembling the edges of $P^1\cup P^2$ it is possible to
replace in $\overline{\cB}$ these two paths with $P^3=[b,a,x,c], 
P^4=[b,c,d,e]$ if $b\neq e$ or $P^5=[a,x,c,b], P^6=[c,d,b,a]$
if $b=e$. Thus, we have again a $(K_v,P_4)$-design.
By the assumptions
on $a$ and $b$ all the vertices of this new design have degree $2$
in at least one block.
\end{proof}

Arguing exactly as in the proof of Theorem \ref{thm:cycle} it is possible to
prove the following result.
\begin{thm}
Let $k>4$. For any admissible $v>1$,
\[\left\{n\geq v-\left\lfloor\frac{k-6}{3}\right\rfloor\ |\ n\equiv0,1\ 
  (\mod{4})\right\}\subseteq
\cL_2P_{k}(v)\subseteq \cL_1P_{k}(v). \]
\end{thm}

\subsection{A construction}
\label{ss:app}
The aim of the current subsection is to provide  for any
admissible $v\equiv 1,2\pmod{4}$ a $(K_v,P_4)$-design $\cB$ with
a vertex having degree $1$ in every block in which it appears.
It  will be then possible to provide a down-link from $\cB$ into
a $(K_{v-1},P_3)$-design $\cB'$, as needed in Theorem \ref{t:81}.
Recall that if $(v-1)(v-2)\not\equiv 0\pmod{4}$,  no $(K_{v-1},P_3)$-design
exists. Thus this condition is necessary for
the existence of a down-link with the required property. We shall 
prove its sufficiency by providing explicit constructions
for all $v\equiv{1,6,9,10}\pmod{12}$.
The approach outlined in Section
\ref{all} shall be extensively used,
by  constructing a partition of the vertices of the
graph $K_v$ in such a way that all the induced complete and complete bipartite
graphs can be down-linked 
to decompositions in $P_3$'s of suitable subgraphs of $K_{v-1}$; these, in turn,
shall yield a decomposition of $\cB'$ with an associated down-link.
\par
Write
$V(K_v)=X_{\ell}\cupdot A_1\cupdot\dots \cupdot A_t$,
where $X_\ell=\{0\}\cup\{1,\dots,\ell-1\}$ for $\ell=6,9,10,13$ and
$|A_i|=12$ for all $i=1,\ldots, t$.
We first construct  a 
$(K_{X_\ell},P_4)$-design $\cB$ which can be down-linked to
a $(K_{X_{\ell}\setminus\{0\}},P_3)$-design $\cB'$. The possible cases are
as follows.
\begin{list}{\labelitemi}{
\leftmargin=0pt}
\item $\ell=6$:
{\scriptsize\begin{align*}
\cB= \{[0,1,2,4],[0,2,3,5],[0,3,4,1],[0,4,5,2],[0,5,1,3]\}
\end{align*}
\begin{align*}
\cB'=\{[1,2,4],[2,3,5],[3,4,1],[4,5,2],[5,1,3]\}
\end{align*}}
\item  $\ell=9$:
  {\scriptsize\begin{align*}
\cB=  \{ & [0,1,2,4],[0,2,3,5],[0,3,4,6],[0,4,5,7],[0,5,6,8],[0,6,7,1], \\ &
[0,7,8,2],[0,8,1,3],[5,8,4,1],[2,5,1,6],[3,6,2,7],[4,7,3,8]\}
\end{align*}
\begin{align*}
  \cB'=\{ &
  [1,2,4],[2,3,5],[3,4,6],[4,5,7],[5,6,8],[6,7,1],[7,8,2],[8,1,3],[8,4,1], \\ &
  [5,1,6],[3,6,2],[7,3,8]\}\cup \{ [2,5,8],[2,7,4]\}
\end{align*}}
\item $\ell=10$:
{\scriptsize\begin{align*}
    \cB=  \{ &
    [0,1,2,4],[0,2,3,5],[0,3,4,6],[0,4,5,7],[0,5,6,8],[0,6,7,9],
    [0,7,8,1],[0,8,9,2],  \\ &
    [0,9,1,3],[1,4,8,2],[2,6,9,4],[4,7,2,5],[5,9,3,7],[7,1,5,8],[8,3,6,1]\}
\end{align*}
{\scriptsize
\begin{align*}
\cB'=\{ &
[1,2,4],[2,3,5],[3,4,6],[4,5,7],[5,6,8],[6,7,9],[7,8,1],[8,9,2],[9,1,3], \\ &
[1,4,8],[6,9,4],[4,7,2],[9,3,7],[7,1,5],[3,6,1] \}\cup
{\scriptsize\{[8,2,6],[2,5,9],[5,8,3]\}}
\end{align*}}}
\item $\ell=13$:
\begin{scriptsize}
\begin{align*}
\cB=  \{ &
[0,1,2,4],[0,2,3,5],[0,3,4,6],[0,4,5,7],[0,5,6,8],[0,6,7,9],[0,7,8,10],
[0,8,9,11], \\ & [0,9,10,12],
[0,10,11,1],[0,11,12,2],[0,12,1,3],[1,4,9,5],[2,5,10,6],[3,6,11,7],\\ &
[4,7,12,8],
[5,8,1,9],
[6,9,2,10],[5,11,3,10],[10,7,1,5],[5,12,9,3],[3,7,2,11],\\ &
[6,12,4,11],[11,8,2,6],
[6,1,10,4],[4,8,3,12]\}.
\end{align*}
\begin{align*}
\cB'=\{ &
[1,2,4],[2,3,5],[3,4,6],[4,5,7],[5,6,8],[6,7,9],[7,8,10],[8,9,11],
[9,10,12],\\ & [10,11,1], [11,12,2],
[12,1,3],[1,4,9],[5,10,6],[3,6,11],[7,12,8],[5,8,1],
[9,2,10],\\ & [5,11,3], [7,1,5],
[5,12,9],[7,2,11], 
[6,12,4],[8,2,6],[6,1,10],[8,3,12]\}\cup\{[9,5,2],\\ &
[11,7,4], [1,9,6],[3,10,7],
[9,3,7], [4,11,8],[10,4,8]\}.
\end{align*}
\end{scriptsize}
\end{list}
We now consider down-links between designs on complete bipartite graphs.
For $X=\{0,1,2\}$ and $Y=\{a,b,c,d,e,f\}$, 
there is a metamorphosis of
the $(K_{X,Y},P_4)$-design
\begin{small}
{\begin{align*}
    \cB=\{[0,a,1,d],[0,b,1,e],[0,c,1,f],[0,d,2,a],[0,e,2,b],[0,f,2,c]\}
\end{align*}}
\end{small}
to the $(K_{X,Y},P_3)$-design
\begin{small}
{\begin{align*}
\cB'=\{[a,1,d],[b,1,e],[c,1,f],[d,2,a],[e,2,b],[f,2,c],
[a,0,b],[c,0,d],[e,0,f]\}.
\end{align*}}
\end{small}
Note that if we remove the paths  $[a,0,b]$,$[c,0,d]$,$[e,0,f]$
from ${\cB'}$, we obtain a bijective down-link  from
$\cB$ to the $(K_{X\setminus\{0\},Y},P_3)$-design
\begin{small}
{\begin{align*}
\cB''=\{[a,1,d],[b,1,e],[c,1,f],[d,2,a],[e,2,b],[f,2,c]\}.
\end{align*}}
\end{small}
Thus, we have actually obtained a
metamorphosis
$\mu:(K_{3,6},P_4)\mbox{-design}\rightarrow(K_{3,6},P_3)\mbox{-design}$ and
a down-link
$\delta:(K_{3,6},P_4)\mbox{-design} \rightarrow (K_{2,6},P_3)\mbox{-design}$.
By gluing together copies of $\mu$ we get metamorphoses
of $P_4$-decompositions into $P_3$-decompositions of
$K_{6,6}$, $K_{9,6}$, $K_{6,12}$, $K_{9,12}$, $K_{12,12}$.
Likewise, using $\delta$ we also determine
down-links  from $P_4$-decompositions of $K_{6,6}$, $K_{6,12}$
and $K_{9,12}$ to  $P_3$-decompositions
of respectively $K_{5,6}$, $K_{5,12}$ and $K_{8,12}$.
For our construction, it will also be necessary to provide 
a metamorphosis of a $(K_{12},P_4)$-design $\cB$ into a
$(K_{12},P_3)$-design $\cB'$.
This is given by
\begin{scriptsize}
\begin{align*}
\cB=\{ & [1,2,3,5],[1,3,4,6],[1,4,5,7],[1,5,6,8],
[1,6,7,9], 
[1,7,8,{10}],  [1,8,9,{11}],[1,9,{10},{12}], \\ & 
[1,{10},{11},{2}],[1,{11},{12},{3}], 
[1,{12},{2},{4}],  [2,{5},{10},{6}],
[3,{6},{11},{7}], [4,{7},{12},{8}], [5,{8},{2},{9}], \\ &
[6,{9},{3},{10}], 
[7,{10},{4},{11}],
[8,{11},{5},{12}],[9,{12},{6},{2}],[{10},{2},{7},{3}], 
[{11},{3},{8},{4}],[{12},{4},{9},{5}]\};
\end{align*}
\begin{align*}
\cB'=\{ & [2,3,5],[3,4,6],[4,5,7],[5,6,8],[6,7,9],[7,8,{10}], 
[8,9,{11}],[9,{10},{12}],
 [{10},{11},{2}], \\ &
[{11},{12},{3}], [{12},{2},{4}], 
[{5},{10},{6}],[{6},{11},{7}], [{7},{12},{8}],
[{8},{2},{9}], [{9},{3},{10}],
[{10},{4},{11}],[{11},{5},{12}],\\ & [{12},{6},{2}],
[{2},{7},{3}],[{3},{8},{4}],
[{4},{9},{5}]\} \cup\{[1,2,5],[1,3,6],[1,4,7],[1,5,8],[1,6,9],
[1,{7},{10}], \\ & 
 [1,{8},{11}], [1,{9},{12}],  [1,{10},{2}],[1,{11},{3}],[1,{12},{4}]\}.
\end{align*}
\end{scriptsize}
Consider now a $(K_v,P_4)$-design with $v=\ell+12t$
where $\ell=1,6,9,10$ and $t>1$.
\begin{list}{\labelitemi}{
\leftmargin=0pt}
\item For $v=1+12t$,
write
\[ K_{1+12t}=(tK_{1,12}\cup tK_{12})\cup\binom{t}{2}K_{12,12}=tK_{13}\cup\binom{t}{2}K_{12,12}. \]
The down-link here is obtained by gluing  down-links from
$P_4$-decompositions
of $K_{13}$ to $P_3$-decompositions of $K_{12}$ with  metamorphoses of
$P_4$-decompositions of  $K_{12,12}$ into  $P_3$-decompositions.
\item For $v=6+12t$, consider
$$K_{6+12t}=K_6\cup tK_{12}\cup tK_{6,12}\cup\binom{t}{2}K_{12,12}.$$
Down-link the $P_4$-decompositions of $K_6$ and $K_{6,12}$ to
respectively $P_3$-decompositions of $K_5$ and $K_{5,12}$ and
consider  metamorphoses of the $P_4$-decompositions of $K_{12}$ and $K_{12,12}$
into $P_3$-decompositions.
\item For $v=9+12t$, let
$$K_{9+12t}=K_9\cup tK_{12}\cup tK_{9,12}\cup\binom{t}{2}K_{12,12}.$$
We know how to down-link the $P_4$-decompositions of $K_{9}$ and $K_{9,12}$ to
$P_3$-decompositions of respectively $K_{8}$
and $K_{8,12}$. As before, there are metamorphoses
of the $P_4$-decompositions  of both $K_{12}$ and $K_{12,12}$ into
$P_3$-decompositions.
\item For $v=10+12t$, observe that
\begin{multline*}
K_{10+12t}=K_{10}\cup tK_{12}\cup tK_{10,12}\cup\binom{t}{2}K_{12,12}=\\
K_{10}\cup tK_{12}\cup tK_{1,12}\cup
tK_{9,12}\cup \binom{t}{2}K_{12,12}=K_{10}\cup tK_{13}\cup tK_{9,12}\cup \binom{t}{2}K_{12,12}.
\end{multline*}
We know how to down-link $P_4$-decompositions of $K_{10}$, $K_{13}$
and $K_{9,12}$ to $P_3$-decompositions  of respectively
$K_{9}$, $K_{12}$ and $K_{8,12}$.
As for the $K_{12,12}$ we argue as in the preceding cases.
\end{list}

\end{document}